\documentclass[11pt,leqno]{article}
\usepackage{amsmath,amssymb}

\begin{document}

\newtheorem{cla}{Claim}[section]
\newtheorem{defi}{Definition}[section]
\newtheorem{rema}{Remark}[section]
\newtheorem{prop}{Proposition}[section]
\newtheorem{lem}{Lemma}[section]
\newtheorem{theo}{Theorem}[section]
\newtheorem{cor}{Corollary}[section]
\newtheorem{conc}{Conclusion}[section]

\author{Pierre GERMAIN}

\title{Finite energy scattering for the Lorentz-Maxwell equation}

\maketitle

\maketitle

\begin{abstract}
In the case where the charge of the particle is small compared to its mass, we describe the asymptotics of the Lorentz-Maxwell equation (Abraham model) for any finite-energy data. As time goes to infinity, we prove that the speed of the particle converges to a certain limit, whereas the electromagnetic field can be decomposed into a soliton plus a free solution of the Maxwell equation.

It is the first instance of a scattering result for general finite energy data in a field-particle equation.
\end{abstract}

\section{Introduction}

\subsection{Presentation of the equation}

\bigskip

\underline{The Abraham model}

\bigskip

The Abraham model describes the interaction of a charged particle with the field that it generates. It consists of a coupling of the Lorentz equation (which governs the movement of the particle), and the Maxwell equation (which gives the evolution of the field).

In the Abraham model, the particle has a fixed, spherically symmetric, charge distribution. This feature is not relativistically invariant, but we will use the relativistic version of the Lorentz equation.

We provide the system with initial data and consider the following Cauchy problem.
\begin{subequations}
\begin{align}
&  m \dot{p}(t) = e E^\rho (q(t)) + e \dot{q}(t) \times B^\rho (q(t)) \;\;\;\;\mbox{with} \;\;\;\; p = \frac{\dot{q}}{\sqrt{1-\dot{q}^2}} \label{AM1} \\
&  \dot{E}(x,t) = \operatorname{curl} B - e \rho(x - q(t)) \dot{q}(t) \label{AM2}\\
&  \dot{B}(x,t) = - \operatorname{curl} E(x,t)  \label{AM3}\\ 
&  \operatorname{div} B(x,t) = 0 \label{AM4}\\
& \operatorname{div} E(x,t) = e \rho(x - q(t)) \label{AM5}\\
&  (q,\dot{q},E,B)_{|t=0} = (q_0,\dot{q}_0,E_0,B_0) \label{AM6}\,\,.
\end{align}
\end{subequations}
The above system is set in the three-dimensional space: $(t,x) \in \mathbb{R} \times \mathbb{R}^3$.
We denote $q$ the position of the particle, $p$ its momentum, and $\rho$ its charge distribution, which is such that
$$
\rho \in \mathcal{C}_0^\infty \;\;\;,\;\;\; \rho \geq 0 \;\;\;,\;\;\;  \rho = \rho(|x|) \;\;\;\,\;\;\;\int_{\mathbb{R}^3} \rho = 1 \,\,.
$$
We also denote, for a function $f$
$$
f^\rho = f * \rho \,\,.
$$
The total charge of the particle is given by $e$, its mass by $m$, and we take all other physical constants (including the speed of light) to be one. The electric and magnetic field are denoted $E$ and $B$.

The initial data are assumed to satisfy the constraints
$$
\operatorname{div} B_0 (x) = 0 \;\;\;\;\;\; \operatorname{div} E_0(x) = e \rho(x-q_0) \,\,.
$$
We can then forget about the conditions on $\operatorname{div} E$ and $\operatorname{div} B$ (equations~(\ref{AM4}) and~(\ref{AM5})) in the above system : if they hold true for the data, this is propagated by the flow given by equations~(\ref{AM1})-(\ref{AM3}).

\bigskip

This model was introduced by Abraham~\cite{abraham03} in 1903 in order to describe the dynamics of the electron. It has recently been intensively studied by mathematicians. Our general reference will be the textbook of Spohn~\cite{spohn04}, in which the interested reader can find further bibliographical indications.

We will be concerned about the asymptotics of the Abraham model. Another important problem, that will not be considered here, is to understand the point particle limit, that is the limit when the support of $\rho$ shrinks to a point.

\bigskip

\underline{Conserved quantities}

\bigskip

The Abraham model actually derives from a Lagrangian, see~\cite{spohn04}. The conservation of the Hamiltonian can also be expressed as the conservation of
$$
\mathcal{E}(\dot{q},E,B) = \frac{m}{\sqrt{1-\dot{q}^2}} + \frac{1}{2}\int_{\mathbb{R}^3} \left( E^2 + B^2 \right) \,\,.
$$
Another conserved quantity is the generalized momentum
$$
\Pi = m p + \int_{\mathbb{R}^3} E \times B \,\,.
$$

\underline{Scaling invariance}

\bigskip

The system~(\ref{AM1})-(\ref{AM6}) is invariant by the following scaling transformation
\begin{equation*}
\begin{aligned}
& E \;\; \longrightarrow \;\; \lambda^{3/2} E(\lambda \cdot,\lambda \cdot) \\
& B \;\; \longrightarrow \;\;\lambda^{3/2} B(\lambda \cdot,\lambda \cdot) \\
& q \;\; \longrightarrow \;\; \lambda^{-1} q(\lambda \cdot) \\
& \dot{q} \;\; \longrightarrow \;\; \dot{q} (\lambda \cdot) \\
& \rho \;\; \longrightarrow \;\; \lambda^{5/2} \rho(\lambda \cdot) \,\,.
\end{aligned}
\end{equation*}
In particular, $\dot{q}$, and the $L^2$ norm of $E$ and $B$, hence the energy $\mathcal{E}$, are left invariant by 
this scaling. In other words, the energy is at the scaling of the equation, the Lorentz-Maxwell system is therefore critical.

\bigskip

\underline{Solitons}

\bigskip

We will call solitons the solutions of the system (\ref{AM1})-(\ref{AM5}) which travel at constant speed, that is
of the form
$$
q = vt \;\;\;\;\; E = e E_v(x-vt) \;\;\;\;\; B = e B_v(x-vt) \,\,,
$$
Thus, the fields $E_v$, $B_v$ solve the elliptic problem
\begin{equation}
\label{eqsoliton}
\left\{
\begin{array}{l}
E_v^\rho (0) + v \times B^\rho_v (0) = 0\\
- v \cdot \nabla E_v = \operatorname{curl} B_v - \rho v \\
-v \cdot \nabla B_v = - \operatorname{curl} E_v\,\,.
\end{array}
\right.
\end{equation}
For any speed $v$ less than 1 in norm there exists only one solution up to space translation (see Spohn~\cite{spohn04}), and it is given by
\begin{equation}
\label{evbv}
E_v(x) = -\nabla \phi_v(x) + v\left(v \cdot \nabla \phi_v(x)\right) \;\;\;\;\;\;\;\; B_v(x) = - v \times 
\nabla \phi_v(x) \,\,,
\end{equation}
where
\begin{equation}
\label{phiv}
\phi_v(x) = \rho(x) * \frac{1}{4\pi \sqrt{(1-v^2)x^2 + (v\cdot x)^2}} \,\,.
\end{equation}

\bigskip

\underline{Related models for field-particle interaction}

\bigskip

We would like to review here some of the models considered in the literature.

It is physically more relevant to allow for a spin of the particle, this possibility is considered in~\cite{imaikinkomechspohn04}.

One can modify the nature of the field interacting with the particle: instead of the electromagnetic field, one can consider a scalar field~\cite{komechspohn98}, a Schr\"odinger field~\cite{komechkopylova06} or a Klein-Gordon field~\cite{imaikinkomechvainberg06}. The qualitative behaviour of all these systems seems to remain essentially the same.

Another possibility is to consider the non-relativistic Lorentz equation~\cite{bambusigalgani93}, 
that is to take as a new definition of the momentum in~(\ref{AM1}) simply $p = \dot{q}$. However, the Maxwell equation remains relativistic and as a consequence there are no solitons that propagate at a speed $v$ larger than $1$.

One can also consider a fully relativistic model by letting the shape of the particle be relativistically invariant: this is the Lorentz model. It is presented in~\cite{spohn04}, but not so many results seem to be available about it.

\subsection{Large time behavior of solutions}

We would like to discuss here known results on large time behavior of solutions.

\bigskip

\underline{Orbital stability and scattering}

\bigskip

It is shown in~\cite{imaikinkomechmauser04} that the Lorentz-Maxwell system exhibits orbital stability, that is to say if data are close to a soliton, the solution remains close to it.

This relies on the following property: for a given generalized momentum $\Pi$, there exists $v$ such that 
the energy is minimized by the solutions $(x_0+vt,v,eE_v(t),eB_v(t))$.

Once orbital stability is proved, the next question is: as $t \rightarrow \infty$, is it possible to decompose the solution into a free electromagnetic field plus a soliton? More specifically this would mean:  there exists $x_0,v \in \mathbb{R}$ and $E_L,B_L$ free solutions of the Maxwell equation such that
\begin{equation}
\label{scattering}
\mbox{as $t \rightarrow \infty$,}\;\;\;\;\;\;\;
\left\{
\begin{array}{l}
\dot{q} \longrightarrow v \\
E \longrightarrow E_v(x-q(t)) + E_L(t) \\
B \longrightarrow B_v(x-q(t)) + B_L(t) 
\end{array}
\right.
\end{equation}
(we leave for the moment unprecise the meaning of the above convergences).

If~(\ref{scattering}) holds, we say that the solution scatters. 

\bigskip

We will see in the following that the solutions scatters under appropriate conditions. This means that the Lorentz-Maxwell equation exhibits the behaviour which is expected for general dispersive field equations: for large time, the solution can be decomposed into solitons and a free solution, that move away from one another. Proving such a large time behaviour is very difficult though for nonlinear field equations. The Lorentz-Maxwell equation is much easier to study, since it involves no nonlinear interaction of the electromagnetic field with itself.
This is one of the great mathematical interests of the Lorentz-Maxwell equation: it provides a tractable model which reproduces some of the features of much more complex situations.

\bigskip

We have seen that orbital stability holds regardless of the ``constants'' of the problem, that is $m$, $e$ and $\rho$.
Let us now discuss under which conditions the solutions scatters.

\bigskip

\underline{Small charge}

\bigskip

If the quotient $\displaystyle \frac{e^2}{m}$ is small, and if the initial fields satisfy
\begin{equation}
\label{e0b0}
|B_0(x)|\,,\,|E_0(x)| \leq \frac{C}{|x|^{3/2+\epsilon}} \;\;\;\mbox{and}\;\;\; |\nabla B_0(x)|\,,\,|\nabla E_0(x)| \leq \frac{C}{|x|^{5/2+\epsilon}} \,\,,
\end{equation}
it is shown in~\cite{spohn04} that scattering occurs. The condition on the moduli of $B_0$, $E_0$ almost includes $L^2$ fields, but the condition on their gradients is very strong.

Our aim in this paper will be to remove these two conditions and prove a scattering result for fields that are merely of finite energy, that is $L^2$.

\bigskip

\underline{Wiener condition}

\bigskip

Another hypothesis under which the solution scatters is the so-called Wiener condition, which requires that
$$
\forall \xi\,\,\,,\,\,\,\widehat{\rho}(\xi) \neq 0 \,\,.
$$
Scattering under this Wiener condition and condition~(\ref{e0b0}) is proved in~\cite{imaikinkomechmauser04}. The beautiful idea, introduced in~\cite{komechspohn97}, is to make use of the following physical fact: a particle which accelerates radiates energy to infinity. Since this amount of radiated energy is bounded, we get a bound on $\ddot{q}$ that enables one to conclude.

\bigskip

\underline{Periodic solutions?}

\bigskip

It has been speculated since the introduction of the Abraham model that periodic in time solutions exist, apparently without any rigorous proof till now, see~\cite{spohn04}. The conditions that we have reviewed under which the solution scatters show that this can happen only for a large charge and with a $\rho$ that violates the Wiener condition.

\section{Main result}

\subsection{Statement}

\begin{theo} 
\label{AMscatters}
If $\displaystyle \frac{e^2}{m}$ is small enough, then for any data of finite energy, that is
$$
\mathcal{E}(\dot{q}_0,E_0,B_0) < \infty\,\,,
$$
the solution of~(\ref{AM1})-(\ref{AM6}) scatters. More precisely, there exists $v_\infty \in \mathbb{R}^3$ (of norm less than 1) and $E_L$, $B_L$ finite-energy solutions of the free Maxwell equation such that
$$
\mbox{as $t \rightarrow \infty$}\;\;,\;\;\;\;\;
\left\{
\begin{array}{l}
\dot{q}(t) \longrightarrow v_\infty \\
E(t) - E_{\dot{q}(t)}(\cdot-q(t)) - E_L(t) \overset{L^2}{\longrightarrow} 0 \\
B(t) - B_{\dot{q}(t)}(\cdot-q(t)) - B_L(t) \overset{L^2}{\longrightarrow} 0 \\
\end{array} \right.
\;\;.
$$
\end{theo}

To our knowledge, this theorem is the first instance of a scattering result for a field-particle interaction equation with general finite energy data; it is also the first instance of such a result where global convergence of the fields is proved.

\subsection{Reduction of the problem}

Consider data as in the above theorem. Using conservation of energy, the global existence of a solution is not hard to prove (we refer to~\cite{spohn04}); conservation of energy also implies that the speed of the particle is bounded away from 1:
\begin{equation}
\label{boundedspeed}
\dot{q}(t) \leq 1 -\epsilon \;\;\;\;,\;\;\;\; \epsilon > 0 \,\,.
\end{equation}

{\bf From now on, we fix initial data $(q_0,\dot{q}_0,E_0,B_0)$ of finite energy. The associated solution $(q,\dot{q},E,B)$ satisfies~(\ref{boundedspeed}).}

\bigskip

Our strategy will be the following: we reduce matters to a linear equation in $\dot{p}$, whose coefficients however depend on $q$, $\dot{q}$. Solving this equation, we obtain the convergence of $\dot{q}$ (or $p$); convergence of the fields follows.

\bigskip

The first step is to introduce the modified fields (following Spohn~\cite{spohn04})
\begin{equation*}
\begin{aligned}
& \bar{E}(x,t) = E(x,t) - e E_{\dot{q}(t)}(x-q(t)) \\
& \bar{B}(x,t) = B(x,t) - e B_{\dot{q}(t)}(x-q(t)) \,\,.
\end{aligned}
\end{equation*}
These new fields somehow measure the distance from the solution to the soliton. In these new coordinates, the 
system~(\ref{AM1})-(\ref{AM6}) becomes
\begin{subequations}
\begin{align}
&  m \dot{p}(t) = e \bar{E}^\rho (q(t)) + e \dot{q}(t) \times \bar{B}^\rho (q(t)) \;\;\;\;\mbox{with} \;\;\;\; p = \frac{\dot{q}}{\sqrt{1-\dot{q}^2}}  \label{MAM1} \\
&  \dot{\bar{E}}(x,t) = \operatorname{curl} \bar{B} (x,t) - e \ddot{q}(t) \nabla_v E_{\dot{q}(t)}(x-q(t)) \label{MAM2} \\
&  \dot{\bar{B}}(x,t) = - \operatorname{curl} \bar{E}(x,t) - e \ddot{q}(t) \nabla_v B_{\dot{q}(t)}(x-q(t)) \label{MAM3}\\ 
& (q,\dot{q},\bar{E},\bar{B})_{|t=0} = (q_0,\dot{q}_0,\bar{E}_0,\bar{B}_0) \label{MAM4} \;\;.
\end{align}
\end{subequations}
Let us now denote $U(t) = \left( \begin{array}{l} U_E \\ U_B \end{array} \right)$ the semi group generated by $\displaystyle \partial_t - \left( \begin{array}{ll} 0 & \operatorname{curl} \\ - \operatorname{curl} & 0 \\ \end{array} \right)$, the equations~(\ref{MAM2})-(\ref{MAM3}) can be rewritten in an integral form as
\begin{equation}
\label{intform}
\begin{aligned}
\left( \begin{array}{l} \bar{E}(t) \\ \bar{B}(t) \end{array} \right) & = U(t) \left( \begin{array}{l} \bar{E}_0 \\ \bar{B}_0 \end{array} \right) + e \int_0^t U(t-s) \ddot{q}(s) \left( \begin{array}{l}  \nabla_v E_{\dot{q}(s)}(x-q(s)) \\ \nabla_v B_{\dot{q}(s)}(x-q(s)) \end{array} \right)\,ds \\
& = \left( \begin{array}{l} \bar{E}_L(t) \\ \bar{B}_L(t) \end{array} \right)
+ e \int_0^t U(t-s) \ddot{q}(s) S(s,x-q(s)) ds \,\,,
\end{aligned}
\end{equation}
setting
$$
\left( \begin{array}{l} \bar{E}_L(t) \\ \bar{B}_L(t) \end{array} \right) = U(t) \left( \begin{array}{l} \bar{E}_0 \\ \bar{B}_0 \end{array} \right) \;\;\;\;\;\;\;\;\;\; S(s,x) = \left( \begin{array}{l}  \nabla_v E_{\dot{q}(s)}(x) \\ \nabla_v B_{\dot{q}(s)}(x) \end{array} \right) \,\,.
$$

So $\displaystyle \left( \begin{array}{l} \bar{E}_L(t) \\ \bar{B}_L(t) \end{array} \right)$ is simply a solution of the free Maxwell equation, whereas $S(s)$ only consists of derivatives of the solitons.

Inserting the above equality in~(\ref{MAM1}), we obtain
\begin{equation*}
\begin{split}
m \dot{p}(t) = e \bar{E}_L^\rho& (q(t)) + e \dot{q}(t) \times \bar{B}_L^\rho(q(t)) \\
& + e^2 \int_0^t  \left[ U_E(t-s) \ddot{q}(s)S^\rho(s) \right](q(t)-q(s))\,ds \\
& + e^2 \dot{q}(t) \times \int_0^t  \left[ U_B(t-s) \ddot{q}(s)S^\rho(s)\right](q(t-q(s))) \,ds \,\,.
\end{split}
\end{equation*}
Finally, the equalities
$$
\dot{q} = \frac{p}{\sqrt{1+p^2}} \overset{\mbox{def}}{=} F(p) \;\;\;\;\;\;\;\;\ddot{q} = F'(p)\dot{p}
$$
give
\begin{equation*}
\begin{aligned}
m \dot{p}(t) = e \bar{E}_L^\rho& (q(t)) + e \dot{q}(t) \times \bar{B}_L^\rho(q(t)) \\
& + e^2 \int_0^t  \left[ U_E(t-s) F'(p(s))  \dot{p}(s) S^\rho(s) \right](q(t)-q(s)) \,ds \\
& + e^2 \dot{q}(t) \times \int_0^t  \left[ U_B(t-s) F'(p(s)) \dot{p}(s) S^\rho(s) \right](q(t)-q(s))\,ds \,\,.
\end{aligned}
\end{equation*}
If one considers any function but $\dot{p}$ as given, we have thus reduced matters to a system which is linear. This will be even more clear setting
\begin{equation}
\label{defabA}
\begin{aligned}
& \alpha(t) = e \bar{E}_L^\rho(q(t)) \\
& \beta(t) = e \bar{B}_L^\rho(q(t)) \\
& A \; : \; f \mapsto \int_0^t  \left[ U_E(t-s) F'(p(s)) f(s) S^\rho(s)\right](q(t)-q(s))\,ds \\
& \;\;\;\;\;\;\;\;\;\;\;\;\;\;+ \dot{q}(t) \times \int_0^t  \left[ U_E(t-s) F'(p(s)) f(s) S^\rho(s) \right](q(t)-q(s))\,ds \,\,.
\end{aligned}
\end{equation}
so that our problem becomes
\begin{equation}
\label{eqalphabeta}
m \dot{p}(t) = \alpha(t) + \dot{q}(t) \times \beta(t) + e^2 (A\dot{p})(t) \,\,.
\end{equation}
We consider $\alpha$, $\beta$, $A$ and $\dot{q}$ as 'fixed' functions (even though they actually depend on the solution of~(\ref{AM1})-(\ref{AM6})); as stressed above, $\dot{p}$ solves a linear equation.

\subsection{Outline of the proof}

We give here the outline of the proof, which will essentially consist on estimates on $\alpha$, $\beta$ and $A$. The actual proof of Propositions~\ref{alphaL2}-\ref{propfields} will be found in Sections~\ref{secalpha}-\ref{secfields}.

As a first and easy step, we shall prove that

\begin{prop}
\label{alphaL2}
If~(\ref{boundedspeed}) holds, the functions $\alpha(t)$ and $\beta(t)$ belong to $L^2([0,\infty))$.
\end{prop}

\begin{prop}
\label{AL2}
If~(\ref{boundedspeed}) holds, the operator $A$ is bounded on $L^2([0,\infty))$.
\end{prop}

If $\frac{e^2}{m}$ is small enough, we thus get
$$
\dot{p} = (m \operatorname{Id} - e^2 A)^{-1} (\alpha + \dot{q} \times \beta) \;\;\;\;\in L^2([0,\infty)) \,\,.
$$

Notice that this gives at once that $\ddot{q} \in L^2([0,\infty))$.
This also implies a better decay of $\alpha$ and $\beta$, as well as better boundedness properties for $A$, as appears in the two next propositions:

\begin{prop}
\label{alphaL1L2}
If~(\ref{boundedspeed}) holds and $\dot{p} \in L^2([0,\infty))$, the functions $\alpha$ and $\beta$ belong to $L^1 ([0,\infty))+ L^2 \cap \partial L^2([0,\infty))$.

(We denote $\partial L^2([0,\infty))$ for functions $f$ such that $f = \dot{g}$, with $g \in L^2$, $g(0)=0$)
\end{prop}

It is easy to see (noticing that $\ddot{q}$ is bounded by the equation) that the last proposition implies that
$$
\alpha + \dot{q} \times \beta \in L^1 + L^2 \cap \partial L^2\,\,.
$$
The next proposition will give us the boundedness of $A$ on this space.

\begin{prop}
\label{AL1L2}
If~(\ref{boundedspeed}) holds and $\dot{p} \in L^2([0,\infty))$, then $A$ is bounded on $L^1([0,\infty)) + L^2 \cap \partial  L^2([0,\infty))$ 
\end{prop}

As a conclusion,
$$
\dot{p} = (m \operatorname{Id} - e^2 A)^{-1} (\alpha + \dot{q} \times \beta) \;\;\;\;\in \left( L^1 + L^2 \cap \partial L^2 \right)([0,\infty)) \,\,.
$$
Now we observe that a function which belongs to $L^2$ together with its derivative has to go to zero. So integrating $\dot{p}$, we get
$$
p(t) \longrightarrow p_\infty \;\;\;\;\mbox{as $t\rightarrow \infty$},
$$
which implies
\begin{equation}
\label{qconverges}
\dot{q}(t) \longrightarrow \dot{q}_\infty \;\;\;\;\mbox{as $t\rightarrow \infty$} \,\,.
\end{equation}
In order to complete the proof of the theorem, it remains to prove that the fields converge as stated.
But since the particle is the only source of the fields, the information that we have gathered on $\dot{p}$ will enable us to conclude.

\begin{prop}
\label{propfields}
If $\dot{p} \in L^1 + L^2 \cap \partial L^2$, there exists a finite energy solution of the free Maxwell equation $(E_L,B_L)$ such that
\begin{equation*}
\mbox{as $t \rightarrow \infty$,}\;\;\;\;\;
\left\{
\begin{array}{l}
E(t) - E_{\dot{q}(t)}(\cdot -q(t)) - E_L(t) \overset{L^2}{\longrightarrow} 0 \\
B(t) - B_{\dot{q}(t)}(\cdot -q(t)) - B_L(t) \overset{L^2}{\longrightarrow} 0 \\
\end{array}
\right.\;\;.
\end{equation*}
\end{prop}

\section{Estimates on $\alpha$ and $\beta$: proofs of Propositions~\ref{alphaL2} and~\ref{alphaL1L2}}

\label{secalpha}

\subsection{The semi group $U(s)$}

Recall the semi-group $U(s)$ is generated by $\displaystyle \partial_t - \left( \begin{array}{ll} 0 & \operatorname{curl} \\ - \operatorname{curl} & 0 \\ \end{array} \right)$. Taking the Fourier transform (that we denote $\widehat{}$\,), it is not hard to see that, denoting
$$
(E_0(t),B_0(t)) = U(t)(E_0,B_0) \,\,,
$$
one has
\begin{equation*}
\begin{aligned}
& \widehat{E_0}(\xi,t) = \cos(t|\xi|) \widehat{E_0} (\xi) + i \sin(t|\xi|) \frac{\xi}{|\xi|} \times \widehat{B}_0 (\xi) \\
& \widehat{B_0}(\xi,t) = \cos(t|\xi|) \widehat{B_0} (\xi) - i \sin(t|\xi|) \frac{\xi}{|\xi|} \times \widehat{E}_0 (\xi) \,\,.
\end{aligned}
\end{equation*}
It is well-known that the inverse Fourier transform of $\frac{\sin(t|\xi|)}{|\xi|}$ (which is the fundamental solution of the scalar wave equation) is given by $\frac{1}{4 \pi t}\delta(|x|-t)$. We deduce from this and the above equality that
\begin{equation}
\label{freeformula}
\begin{split}
& E_0(x,t) = \frac{1}{4 \pi t^2} \int_{|y-x|=t} E_0(y)\,dy + \frac{1}{4 \pi t} \int_{|y-x|=t} \frac{y-x}{|y-x|}\cdot\nabla E_0 \,dy \\
&\;\;\;\;\;\;\;\;\;\;\;\;\;\;\;\;\;\;\;\;\;\;\;\;\;\;\;\;\;\;\;\;\;+ \frac{1}{4 \pi t} \int_{|y-x|=t} \operatorname{curl} B_0(y) \,dy \\
& B_0(x,t) = \frac{1}{4 \pi t^2} \int_{|y-x|=t} B_0(y)\,dy + \frac{1}{4 \pi t} \int_{|y-x|=t} \frac{y-x}{|y-x|} \cdot \nabla B_0 \,dy \\
&\;\;\;\;\;\;\;\;\;\;\;\;\;\;\;\;\;\;\;\;\;\;\;\;\;\;\;\;\;\;\;\;\;- \frac{1}{4 \pi t} \int_{|y-x|=t} \operatorname{curl} E_0(y) \,dy 
\end{split}
\end{equation}
(the above integrals over surfaces are understood with respect to the standard surface measure).

\subsection{A change of variable}

Combining the above formula and~(\ref{defabA}), that is the definition of $\alpha$ and $\beta$, we see that in order to prove Propositions~\ref{alphaL2} and~\ref{alphaL1L2}, we have to study surface integrals of functions which are $L^2$ in space. More specifically, $\alpha$ and $\beta$ can be written as linear combinations of functions of the type
\begin{equation}
\begin{aligned}
& g_1(t) = \frac{1}{t^2} \int_{|z-q(t)|=t} f^\rho (z)\,dz \\
& g_2(t) = \frac{1}{t} \int_{|z-q(t)|=t}  \frac{z-q(t)}{|z-q(t)|} \cdot \nabla f^\rho(z) \,dz \\
& g_3(t) = \frac{1}{t} \int_{|z-q(t)|=t} \operatorname{curl} f^\rho(z)\,dz\,\,,
\end{aligned}
\end{equation}
where $f$ is an $L^2$ function. Instead of proving estimates on $\alpha$ and $\beta$, we will prove estimates on $g_1$ and $g_2$, the case of $g_3$ being very similar.

\bigskip

Since our aim will be to integrate in time, it appears, considering the definition of the $g_i$, that the following change of variable is natural:
\begin{equation}
\label{coc}
\begin{array}{llll}
 \phi : & \mathbb{R}^3 & \rightarrow & \mathbb{R}^3 \\
 & y = s\omega & \mapsto & z = q(s) + s\omega 
\end{array}\,\,,
\end{equation}
where $s \in \mathbb{R}^+$, $\omega \in \mathbb{S}^2$, in other words, $(s,\omega)$ are the polar coordinates of $y$.

Assuming that $\omega = e_1$ (the first vector of an orthonormal basis $(e_1,e_2,e_3))$, we get
\begin{equation}
\label{boundcoc}
\begin{aligned}
& \frac{\partial s}{\partial z^1} = \frac{1}{1+\dot{q}^1} \;\;\;\;\;\; \frac{\partial s}{\partial z^2} = \frac{\partial s}{\partial z^3} = 0 \\
& \frac{\partial \omega^1}{\partial z^i} = 0 \;\;\;\;\;\;\mbox{for any $i$} \\
& \frac{\partial \omega^2}{\partial z^1} = -\frac{1}{s}\frac{\dot{q}^2}{1+\dot{q}^1} \;\;\;\;\;\; \frac{\partial \omega^2}{\partial z^2} = \frac{1}{s} \;\;\;\;\;\; \frac{\partial \omega^3}{\partial z^3} = 0 \\
& \frac{\partial \omega^3}{\partial z^1} = -\frac{1}{s}\frac{\dot{q}^3}{1+\dot{q}^1} \;\;\;\;\;\; \frac{\partial \omega^3}{\partial z^2} = 0 \;\;\;\;\;\; \frac{\partial \omega^3}{\partial z^3} = \frac{1}{s} \\
& \frac{\partial y}{\partial z} = \left( \displaystyle \begin{array}{ccc} \frac{1}{1+\dot{q}^1} & -\frac{\dot{q}^2}{1+\dot{q}^1} & -\frac{\dot{q}^3}{1+\dot{q}^1} \\ 0 & 1 & 0 \\ 0 & 0 & 1 \end{array} \right)
\end{aligned}
\end{equation}
This implies in particular that
$$
\left| \frac{\partial y}{\partial z} \right| = \frac{1}{1+ \langle \dot{q}(s)\,,\,\omega \rangle} \,\,.
$$
It is at this point that~(\ref{boundedspeed}) plays a crucial role: due to this inequality, $\phi$ is a global diffeomorphism.

\subsection{Proof of Proposition~\ref{alphaL2}}

\label{section33}

As we saw above, we reduce matters to studying $g_1$ and $g_2$. 
Thus to prove Proposition~\ref{alphaL2}, it suffices to prove that $t \mapsto g_1(t)$ and $t \mapsto g_2(t)$ belong to $L^2$.

We will prove this for $g_2$ - it being easier for $g_1$. By definition of $g_2$,
\begin{equation*}
\begin{split}
\|g_2\|_{L^2([0,\infty))}^2 & = \int_0^\infty \left( \frac{1}{t} \int_{|z-q(t)|=t}  \frac{z-q(t)}{|z-q(t)|} \cdot \nabla f^\rho(z) \,dz \right)^2 \,dt \\
& \leq C \int_0^\infty \int_{|z-q(t)|=t} |\nabla f^\rho(z)|^2 \,dz \,dt \;\;\;\;\mbox{by H\"older's inequality} \\
& = C \int_{\mathbb{R}^3} |\nabla f^\rho(\phi(y))|^2 \,dy \;\;\;\;\mbox{where $\phi$ is defined in~(\ref{coc})} \\
& = C \int_{\mathbb{R}^3} |\nabla f^\rho(z)|^2 \frac{dz}{1 + \langle \dot{q}(s)\,,\,\omega \rangle} \;\;\;\;\mbox{changing coordinates} \\
& \leq C \|\nabla f^\rho \|_2^2 \;\;\;\;\mbox{using~(\ref{boundedspeed})}\,\,.
\end{split}
\end{equation*}

\subsection{Proof of Proposition~\ref{alphaL1L2}}

As in the previous subsection, we prove estimates on $g_1$ and $g_2$ instead of $\alpha$ and $\beta$. Thus our aim here will be to prove that $g_1$ and $g_2$ belong to $L^1 + \partial L^2$.

\bigskip

The result for $g_1$ is easily obtained: a straightforward estimate as in the previous subsection yields: $g_1 \in L^1$.

\bigskip

Now we would like to prove that $g_2$ belongs to $L^1 + \partial L^2$. The idea is to integrate $g_2$ in time, and try and write the result as the integral of an $L^1$ function, plus an $L^2$ function. Differentiating this equality, we will get the desired result. 

Let us introduce the following notation: $B_t$ is the Euclidean ball of radius $t$, and $n_t$ the exterior normal to $\phi(\partial B_t)$. By definition of $g_2$,
\begin{equation*}
\begin{aligned}
\int_0^T g_2(t) \,dt & = \int_0^T \frac{1}{t} \int_{|z-q(t)|=t}  \frac{z-q(t)}{|z-q(t)|} \cdot \nabla f^\rho(z) \,dz \,dt \\
& = \int_{\phi(B_T)} \frac{y}{|y|^2} \cdot \nabla f^\rho(z) \,\frac{dz}{1+\langle \dot{q}(s)\,,\,\omega\rangle}
\;\;\;\;\;\;\mbox{proceeding as in~\ref{section33}} \\
& = \int_{\phi(\partial B_T)} n_t(z) \cdot \frac{y}{|y|^2} f^\rho(z) \,\frac{dz}{1+\langle \dot{q}(s)\,,\,\omega \rangle} \\
& \;\;\;\;\;\;\;\;- \int_{\phi(B_T)} f^\rho(z) \nabla_z \cdot \left( \frac{y}{|y|^2} \frac{1}{1+\langle \dot{q}(s)\,,\,\omega\rangle} \right) \,dz  \\
& = \int_{\phi(\partial B_T)} n_t(z) \cdot \frac{y}{|y|^2} f^\rho(z) \,\frac{dz}{1+\langle \dot{q}(s)\,,\,\omega \rangle} \\
& \;\;\;\;\;\;\;\;- \int_{B_T} f^\rho(\phi(y)) \cdot \nabla_z \left( \frac{y}{|y|^2} \frac{1}{1+\langle \dot{q}(s)\,,\,\omega\rangle} \right) \left[1+\langle \dot{q}(s)\,,\,\omega \rangle \right]\,dy \\
& \overset{\mbox{def}}{=} K(T) - \int_{B_T} L(y)\,dy \,\,.
\end{aligned}
\end{equation*}
(in the above, $z$, $y$, $s$, and $\omega$ are functions of $y$ or $z$, given by~(\ref{coc})). Differentiating in time,
$$
g_2(t) = \dot{K}(t) - \int_{\partial B_t} L(y) \,dy \,\,.
$$
This is the decomposition of $g_2$ that we were looking for. Indeed, we have the following 

\begin{cla}
With $K$ and $L$ defined above, one has
$$
\dot{K} \in L^2 \cap \partial L^2 \;\;\;\mbox{and}\;\;\;\; t \mapsto \int_{\partial B_t} L \;\;\in L^1 \,\,.
$$
\end{cla}

\textsc{Proof of the claim:} First, using~(\ref{boundcoc}), we see that
\begin{equation*}
\begin{split}
|L(y)| & \leq C |f^\rho(\phi(y))| \left| \nabla_z \left( \frac{y}{|y|^2} \frac{1}{1+\langle \dot{q}(s)\,,\,\omega\rangle} \right) \right| \\
& \leq C |f^\rho(\phi(y))| \left(\frac{1}{|y|^2} + \frac{|\ddot{q}(|y|)|}{|y|} \right)\,\,.
\end{split}
\end{equation*}
But $f^\rho \circ \phi \in L^2\cap L^\infty (\mathbb{R}^3)$, and, since $\ddot{q} \in L^2([0,\infty))$,
$$
y \mapsto \left(\frac{1}{|y^2|} + \frac{|\ddot{q}(|y|)|}{|y|} \right) \;\;\in L^1 + L^2(\mathbb{R}^3)\,\,.
$$
Thus $L \in L^1(\mathbb{R}^3)$, which implies that $t \mapsto \int_{\partial B_t} L \in L^1([0,\infty))$.
\smallskip

The definition of $K$ ensures as in~\ref{section33} that it belongs to $L^2([0,\infty))$, and that $K(0)=0$. 
Finally, $\dot{K}$ belongs to $L^2$ since it can be written as the sum of two $L^2$ functions
$$
\dot{K}(t) = g_2(t) + \int_{\partial B_t} L(y) \,dy \,\,.\;\;\;\;\;\;\;\blacksquare
$$

\section{Estimates on $A$: proofs of Propositions~\ref{AL2} and~\ref{alphaL1L2}}

\subsection{Estimates on the kernel of $A$}

\label{secA}

First, we will need in the following to estimate derivatives of the soliton. Equations~(\ref{evbv}) and~(\ref{phiv}) imply easily
\begin{equation}
\label{derivsoliton}
\left| \nabla_x^k \nabla_v^\ell E_v(x)\right| \;\;,\;\;\left| \nabla_x^k \nabla_v^\ell B_v(x)\right| \leq \frac{C}{|x|^{k+2}} \,\,.
\end{equation}

Recall that $A$ is defined in~(\ref{defabA}). Let us write it as a kernel operator, that is
$$
Af(t) = \int_0^t a(t,s) f(s) \,ds \,\,.
$$
The actual formula for $a$ is long and involved. We do not keep track of details that are not relevant for the analysis and thus consider in the following the somewhat simplified version
$$
a(t,s) = \dot{q}(t)  F'(p(s)) \left[U_E(t-s) S^\rho(s) \right](q(t)-q(s))\,\,.
$$
The following lemma will provide us with bounds on $a$ and its derivatives

\begin{lem}
\label{lema}
Denoting $\langle x \rangle = \sqrt{1+x^2}$, one has
\begin{equation*}
\begin{aligned}
& \left| \left[U_E(t-s) S^\rho(s) \right](q(t)-q(s))\right| \leq \frac{C}{\langle t-s \rangle^2} \\
& \left| \left[U_E(t-s) \nabla_x S^\rho(s) \right](q(t)-q(s))\right| \leq \frac{C}{\langle t-s \rangle^3} \\
\end{aligned}
\end{equation*}
In particular,
$$
a(t,s) \leq \frac{C}{\langle t-s \rangle^2} \,\,.
$$
\end{lem}

\textsc{Proof of the lemma:} 
We prove only the first inequality, the other being similar.
By formula~(\ref{freeformula}), it suffices to prove the desired bound for functions of the type
$$
\frac{1}{(t-s)^2}\int_{|y-(q(t)-q(s))|=t-s} g \;\;\;\;\;\mbox{or}\;\;\;\;\; \frac{1}{t-s} \int_{|y-(q(t)-q(s))|=t-s} \nabla g
$$
where $g = \nabla_v E_v$ or $\nabla_v B_v$. And this follows from~(\ref{boundedspeed})~(\ref{derivsoliton}). $\blacksquare$

\subsection{Proof of Proposition~\ref{AL2}}

We aim at proving the boundedness of $A$ on $L^2([0,\infty))$.

Unfortunately, $a(t,s)$ is not a function of $(t-s)$, so it is not a convolution operator. However,  Lemma~\ref{lema} asserts that is possible to bound $a$ by an integrable function of $(t-s)$.

This gives at once the boundedness of $A$ on $L^2([0,\infty))$.

\subsection{A small lemma on $L^2$ functions}

We will need the following elementary lemma.

\begin{lem}
\label{l2lemma}
If $F' \in L^2([0,\infty))$, the function
$$
s \mapsto \int_0^\infty \frac{|F(\tau+s)-F(s)|}{\langle \tau \rangle^3} d\tau
$$
also belongs to $L^2([0,\infty))$.
\end{lem}

\textsc{Proof:} Integrating the square of the above function gives
\begin{equation}
\begin{aligned}
\int_0^\infty &\left( \int_0^\infty  \frac{|F(\tau+s)-F(s)|}{\langle \tau \rangle^3} d\tau \right)^2 \,ds \\
& \leq C \int_0^\infty \int_0^\infty \frac{|F(\tau+s)-F(s)|^2}{\langle \tau \rangle^4} \,d\tau \,ds \;\;\;\;\;\mbox{by H\"older's inequality}\\
& \leq C \int_0^\infty \int_0^\infty \int_{s}^{\tau+s} |F'(\sigma)|^2 \,d\sigma \frac{d\tau}{\langle \tau \rangle^3} ds \;\;\;\;\;\mbox{also by H\"older's inequality}\\
& = C \int_0^\infty \int_0^\infty \int_0^\tau |F'(\sigma+s)|^2\,d\sigma \frac{ d \tau }{ \langle \tau \rangle^3} \,ds \\
& \leq C \int_0^\infty \int_0^\tau \|F'\|_{L^2}^2 \,d\sigma \frac{d\tau}{\langle \tau \rangle^3} \;\;\;\;\;\mbox{by Fubini's theorem}\\
& \leq C \|F'\|_{L^2}^2 \int_0^\infty \frac{d\tau}{\langle \tau \rangle^2} < \infty \,\,.\;\;\;\;\blacksquare
\end{aligned}
\end{equation}

\subsection{Proof of Proposition~\ref{AL1L2}}

We would like to prove that $A$ is bounded on $L^1([0,\infty)) + L^2 \cap \partial L^2([0,\infty))$. Lemma~\ref{lema} implies at once that $a$ is bounded on $L^1([0,\infty))$.

So let us consider $f \in \partial L^2 \cap L^2 ([0,\infty))$, that is
$$
f=\dot{g}\;\;\;\;\mbox{with}\;\;\;\;g \in L^2\;\;,\;\;\dot{g}\in L^2\;\;,\;\;g(0)=0\,\,.
$$
We observe that, up to a remainder term, we can interchange the derivatives of $a$ with respect to $s$ and $t$
\begin{equation}
\label{ats}
\begin{split}
\frac{\partial a}{\partial s}(t,s) = - \frac{\partial a}{\partial t} & + \dot{q}(t) F'(p(s)) (\dot{q}(t) - \dot{q}(s)) \left[U_E(t-s) \nabla_x S^\rho(s) \right](q(t)-q(s)) \\
& + \dot{q}(t) \dot{p}(s) F''(p(s)) \left[U_E(t-s) S^\rho(s) \right](q(t)-q(s)) \\
& + \dot{q}(t) F'(p(s)) \ddot{q}(s) \left[U_E(t-s) \nabla_v S^\rho(s) \right](q(t)-q(s)) \\
& + \ddot{q}(t) F'(p(s)) \left[U_E(t-s) S^\rho(s) \right](q(t)-q(s)) \,\,, \\
= - \frac{\partial a}{\partial t} &  + R \,\,,
\end{split}
\end{equation}
where, by Lemma~\ref{lema} and since $\ddot{q} \in L^2 \cap L^\infty$, we can bound the remainder by
\begin{equation}
\label{boundR}
|R(t,s)| \leq C \frac{|\dot{q}(t)-\dot{q}(s)|}{ \langle t-s \rangle^3} + \frac{C h(s)}{ \langle t-s \rangle^2} + \frac{C h(t)}{ \langle t-s \rangle^2} \;\;\;\;\;\mbox{with}\;\;h \in L^2\cap L^\infty \,\,.
\end{equation}
Now we perform an integration by parts and use~(\ref{ats}) to get
\begin{equation}
\label{onetwothree}
\begin{split}
Af(t) & = \int_0^t a(t,s) \dot{g}(s)\,ds \\
& = a(t,t)g(t) - a(t,0)g(0) - \int_0^t \frac{\partial a}{\partial s}(t,s) g(s) \,ds \\
& = a(t,t)g(t) - a(t,0)g(0) + \int_0^t \frac{\partial a}{\partial t}(t,s) g(s) \, ds - \int_0^t R(t,s) g(s) \, ds\\
& = - a(t,0)g(0) - \int_0^t R(t,s) g(s) + \frac{\partial}{\partial t} \left( \int_0^t a(t,s) g(s) \, ds \right)  \, ds \\
& = I - II + III\,\,.
\end{split}
\end{equation}
Let us check that the above sum belongs to $L^1 + L^2 \cap \partial L^2$. 

\bigskip

The decay of $a$ ensures first that $I : t \mapsto a(t,0)g(0)$ belongs to $L^1$. 

\bigskip

Let us now consider $II$. By~(\ref{boundR}), it can be bounded by
$$
C \int_0^t \left( \frac{|\dot{q}(t)-\dot{q}(s)|}{ \langle t-s \rangle^3} + \frac{h(s)}{ \langle t-s \rangle^2} + \frac{h(t)}{ \langle t-s \rangle^2} \right) g(s) \,ds \,\,,
$$
We now check that each of the three pieces of the right-hand side belong to $L^1$:
\begin{itemize}
\item We begin with the first one.
\begin{equation}
\begin{split}\;\;\;\;\;\;\;\;
\int_0^\infty \int_0^t & \frac{|\dot{q}(t)-\dot{q}(s)|}{ \langle t-s \rangle^3} g(s)\, ds \,dt \\
&  = \int_0^\infty g(s) \int_0^\infty \frac{\dot{q}(\tau +s) - \dot{q}(s)}{\langle \tau \rangle^3} \, d\tau \,ds
\leq C \|\ddot{q}\|_2 \|g\|_2
\end{split}
\end{equation}
by Lemma~\ref{l2lemma}.
\item The second one can be written as
$$
\frac{1}{\langle s \rangle^2} * (hg)
$$
and thus it belongs to $L^1$.
\item The third one can be written as
$$
h \left( \frac{1}{\langle s \rangle^2} * g\right)
$$
which also belongs to $L^1$.
\end{itemize}

\bigskip

It only remains to deal with $III$. Since $g$ belongs to $L^2$, so does $\displaystyle t \mapsto \int_0^t a(t,s) g(s) \, ds$. As a conclusion, $III$ belongs to $\partial L^2 ([0,\infty))$. Finally, $III$ belongs to $L^2$ since we can write it as a sum of $L^2$ functions:
$$
III = Af(t) - I + II\,\,.
$$

\section{Scattering for the fields: proof of Proposition~\ref{propfields}}
\label{secfields}

First, we observe that since $\dot{p} \in L^1 + L^2 \cap \partial L^2$, there holds
$$
\ddot{q} \in L^1 + \partial L^2 \,\,.
$$
Indeed, we have $\dot{q} = F(p)$, and $\dot{p} = f + \dot{g}$, with $f \in L^1$, $g \in L^2$, $\dot{g} \in L^2$.
Thus
$$
\ddot{q} = (f + \dot{g}) F'(p) = f F'(p) + \left( g F'(p)\right)^{\cdot} - g \dot{p} F''(p) \;\;\;\;\in L^1 + \partial L^2\,\,.
$$
Coming back to~(\ref{intform}), we see that in order to prove Proposition~\ref{propfields}, it suffices to show that
$$
\int_0^\infty U(-s) \ddot{q}(s) \left( \begin{array}{l}  \nabla_v E_{\dot{q}(s)}(x-q(s)) \\ \nabla_v B_{\dot{q}(s)}(x-q(s)) \end{array} \right)\,ds \;\;\in L^2(\mathbb{R}^3) 
$$
and
$$
\mbox{as $t\rightarrow \infty$,}\;\;\;\;\int_t^\infty U(t-s) \ddot{q}(s) \left( \begin{array}{l}  \nabla_v E_{\dot{q}(s)}(x-q(s)) \\ \nabla_v B_{\dot{q}(s)}(x-q(s)) \end{array} \right)\,ds \overset{L^2}{\longrightarrow} 0 \,\,.
$$
We shall focus on the convergence to zero, the other point being very similar. 
Replacing in the above expression $\ddot{q}(s)$ by $f + \dot{g}$, with $f \in L^1$, $g \in L^2 \cap \partial L^2$, we want to prove the convergence to zero of
$$
\int_t^\infty U(t-s) \left(f(s) + \dot{g}(s)\right) \left( \begin{array}{l}  \nabla_v E_{\dot{q}(s)}(x-q(s)) \\ \nabla_v B_{\dot{q}(s)}(x-q(s)) \end{array} \right)\,ds \,\,.
$$
The term involving $f$ is easily dealt with, since $f \in L^1$. As for the term involving $\dot{g}$, integration by parts combined with equation~(\ref{eqsoliton}) yields
\begin{equation*}
\begin{split}
\int_t^\infty&  U(t-s) \dot{g}(s)  \left( \begin{array}{l}  \nabla_v E_{\dot{q}(s)}(x-q(s)) \\ \nabla_v B_{\dot{q}(s)}(x-q(s)) \end{array} \right)\,ds \\
& = - g(t) \left( \begin{array}{l}  \nabla_v E_{\dot{q}(t)}(x-q(t)) \\ \nabla_v B_{\dot{q}(t)}(x-q(t)) \end{array} \right) \\
& \;\;\;\;\;- \int_t^\infty  U(t-s) (g(s) \otimes \ddot{q}(s)) \nabla_v \left( \begin{array}{l}  \nabla_v E_{\dot{q}(s)}(x-q(s)) \\ \nabla_v B_{\dot{q}(s)}(x-q(s)) \end{array} \right)\,ds \\
& \;\;\;\;\;+ \int_t^\infty  U(t-s) g(s) \left( \begin{array}{l} - \nabla_x E_{\dot{q}(s)}(x-q(s)) + \rho(x-q(s)) \operatorname{Id} \\ \nabla_x B_{\dot{q}(s)}(x-q(s)) \end{array} \right) \,ds \,\,.
\end{split}
\end{equation*}
Since $g(t)$ converges to zero and $g(s) \ddot{q}(s) \in L^1$, matters reduce to showing that the third term in the right-hand side goes to zero. We will actually show that
\begin{equation}
\label{l2l2}
F(x) \overset{\mbox{def}}{=} \int_0^\infty  \left| U(-s) g(s) \left( \begin{array}{l}  \nabla_x E_{\dot{q}(s)}(x-q(s)) + \rho(x-q(s)) \operatorname{Id} \\ \nabla_x B_{\dot{q}(s)}(x-q(s)) \end{array} \right) \right| \,ds  \,\,.
\end{equation}
belongs to $L^2(\mathbb{R}^3)$. We can bound uniformly the norm of 
$$\left( \begin{array}{l}  \nabla_x E_{\dot{q}(s)}(x-q(s)) + \rho(x-q(s))\operatorname{Id} \\ \nabla_x B_{\dot{q}(s)}(x-q(s)) \end{array} \right)$$ 
by $\frac{1}{\langle x - q(s) \rangle^3}$, and its derivative by $\frac{1}{\langle x - q(s) \rangle^4}$. After applying $U(s)$, formula~(\ref{freeformula}) gives a bound of the form
\begin{equation}
\begin{split}
& \left| U(-s) \left( \begin{array}{l}  \nabla_x E_{\dot{q}(s)}(x-q(s)) + \rho(x-q(s))\operatorname{Id} \\  \nabla_x B_{\dot{q}(s)}(x-q(s)) \end{array} \right) \right| \\
& \;\;\;\;\;\;\;\;\;\;\;\;\;\;\;\;\;\;\;\;\leq \frac{C}{s^2} \int_{|z|=s} \frac{dz}{\langle x-q(s)- z \rangle^3} + \frac{C}{s} \int_{|z|=s} \frac{dz}{\langle x- q(s)-z \rangle^4} \,\,.
\end{split}
\end{equation}
As usual, the treatments of the two terms on the right-hand side are very similar, so we keep only the second one, which is slightly harder, and estimate
\begin{equation}
\begin{split}
|F(x)|^2 & \leq C \left( \int_0^\infty g(s) \frac{1}{s} \int_{|z|=s} \frac{dz}{\langle x - z -q(s) \rangle^4} \,ds \right)^2 \\
 & = C \left( \int_{\mathbb{R}^3} \frac{g(|z|)}{|z|} \frac{1}{\langle x - z - q(|z|) \rangle^4} \,dz \right)^2 \\
& \leq C \int_{\mathbb{R}^3} \frac{g(|z|)^2}{|z|^2} \frac{1}{\langle x - z - q(|z|) \rangle^4} \,dz \int_{\mathbb{R}^3} \frac{1}{\langle x - z - q(z) \rangle^4} \,dz  \\
& \leq C  \int_{\mathbb{R}^3} \frac{g(|z|)^2}{|z|^2} \frac{1}{\langle x - z - q(|z|) \rangle^4} \,dz \;\;\;\;\mbox{by~(\ref{boundedspeed})} \,\,.
\end{split}
\end{equation}
To conclude, it suffices to integrate in $x$, use Fubini, and observe that $\frac{g(|z|)}{|z|} \in L^2 (\mathbb{R}^3)$:
$$
\int_{\mathbb{R}^3} |F(x)|^2\,dx \leq C \int_{\mathbb{R}^3}  \frac{g(|z|)^2}{|z|^2} \int_{\mathbb{R}^3} \frac{1}{\langle x - z - q(z) \rangle^4} \,dx\,dz <\infty \,\,.
$$

\bigskip

\bigskip

{\bf Acknowledgements:} The author is very grateful to Igor Rodnianski for suggesting this problem to him, and very interesting conversations about it.

\bigskip

Pierre GERMAIN 

{\sc Courant Institute of Mathematical Sciences, New York University, New York, NY 10012-1185, USA}

\bigskip
  
{\tt pgermain@math.nyu.edu}

\end{document}